\font\eusm=eusm10                   


\font\eusms=eusm7                       

\font\eusmss=eusm5                      

\font\scriptsize=cmr7

\input amstex

\documentstyle{amsppt}
  \magnification=1100
  \hsize=6.2truein
  \vsize=9.0truein
  \hoffset 0.1truein
  \parindent=2em

\NoBlackBoxes

\newcount\mycitestyle \mycitestyle=1



\newcount\theTime
\newcount\theHour
\newcount\theMinute
\newcount\theMinuteTens
\newcount\theScratch
\theTime=\number\time
\theHour=\theTime
\divide\theHour by 60
\theScratch=\theHour
\multiply\theScratch by 60
\theMinute=\theTime
\advance\theMinute by -\theScratch
\theMinuteTens=\theMinute
\divide\theMinuteTens by 10
\theScratch=\theMinuteTens
\multiply\theScratch by 10
\advance\theMinute by -\theScratch

\def\today{{\number\day\space
 \ifcase\month\or
  January\or February\or March\or April\or May\or June\or
  July\or August\or September\or October\or November\or December\fi
 \space\number\year}}





\define\ah{{\hat a}}                         














\define\biggnm#1{                            
  \bigg|\bigg|#1\bigg|\bigg|}

\define\bignm#1{                             
  \big|\big|#1\big|\big|}

\define\bh{\hat b}                           










\define\Cpx{\bold C}                         








\define\eh{\hat e}                           



\define\ev{{\text{\rm ev}}}                  




\define\fh{\hat f}                           

\define\fpamalg#1{{\dsize\;                  
     \operatornamewithlimits*_{#1}\;}}


\define\freeprod#1#2{\mathchoice             
     {\operatornamewithlimits{\ast}_{#1}^{#2}}
     {\raise.5ex\hbox{$\dsize\operatornamewithlimits{\ast}
      _{#1}^{#2}$}\,}
     {\text{oops!}}{\text{oops!}}}

\define\freeprodi{\mathchoice                
     {\operatornamewithlimits{\ast}
      _{\iota\in I}}
     {\raise.5ex\hbox{$\dsize\operatornamewithlimits{\ast}
      _{\sssize\iota\in I}$}\,}
     {\text{oops!}}{\text{oops!}}}










\define\GNS{{\text{\rm GNS}}}                



\define\Hil{{\mathchoice                     
     {\text{\eusm H}}
     {\text{\eusm H}}
     {\text{\eusms H}}
     {\text{\eusmss H}}}}







\define\Hilo{{\overset{\scriptsize o}        
     \to\Hil}}







\define\id{\text{\rm id}}                    






\define\KHil{{\mathchoice                    
     {\text{\eusm K}}
     {\text{\eusm K}}
     {\text{\eusms K}}
     {\text{\eusmss K}}}}







\define\lrnm#1{\left|\left|#1\right|\right|} 










\define\nm#1{||#1||}                         


\define\Naturals{{\bold N}}                  




\define\onehat{\hat 1}                       



\define\oup{^{\text{\rm o}}}                 

\define\owedge{{                             
     \operatorname{\raise.5ex\hbox{\text{$
     \ssize{\,\bigcirc\llap{$\ssize\wedge\,$}\,}$}}}}}

\define\owedgeo#1{{                          
     \underset{\raise.5ex\hbox
     {\text{$\ssize#1$}}}\to\owedge}}










\define\QED{\newline                         
            \line{$\hfill$\qed}\enddemo}











\define\restrict{\lower .3ex                 
     \hbox{\text{$|$}}}




\define\smd#1#2{\underset{#2}\to{#1}}          

\define\smdb#1#2{\undersetbrace{#2}\to{#1}}    

\define\smdbp#1#2#3{\overset{#3}\to            
     {\smd{#1}{#2}}}

\define\smdbpb#1#2#3{\oversetbrace{#3}\to      
     {\smdb{#1}{#2}}}

\define\smdp#1#2#3{\overset{#3}\to             
     {\smd{#1}{#2}}}

\define\smdpb#1#2#3{\oversetbrace{#3}\to       
     {\smd{#1}{#2}}}

\define\smp#1#2{\overset{#2}\to                
     {#1}}                                     








\define\Tcirc{\bold T}                       

\define\Teu{\text{\eusm T}}


\define\tocdots                              
  {\leaders\hbox to 1em{\hss.\hss}\hfill}    

\define\tr{\text{\rm tr}}                    




























  \newcount\bibno \bibno=0
  \def\newbib#1{\advance\bibno by 1 \edef#1{\number\bibno}}
  \ifnum\mycitestyle=1 \def\cite#1{{\rm[\bf #1\rm]}} \fi


  \newcount\ignorsec \ignorsec=0
  \def\notasec{\ignorsec=1}

  \newcount\secno \secno=0
  \def\newsec#1{\procno=0 \subsecno=0 \ignorsec=0
    \advance\secno by 1 \edef#1{\number\secno}
    \edef\currentsec{\number\secno}}

  \newcount\subsecno
  \def\newsubsec#1{\procno=0 \advance\subsecno by 1 \edef#1{\number\subsecno}
    \edef\currentsec{\number\secno.\number\subsecno}}

  \newcount\appendixno \appendixno=0
  \def\newappendix#1{\procno=0 \ignorsec=0 \advance\appendixno by 1
    \ifnum\appendixno=1 \edef\appendixalpha{\hbox{A}}
      \else \ifnum\appendixno=2 \edef\appendixalpha{\hbox{B}} \fi
      \else \ifnum\appendixno=3 \edef\appendixalpha{\hbox{C}} \fi
      \else \ifnum\appendixno=4 \edef\appendixalpha{\hbox{D}} \fi
      \else \ifnum\appendixno=5 \edef\appendixalpha{\hbox{E}} \fi
      \else \ifnum\appendixno=6 \edef\appendixalpha{\hbox{F}} \fi
    \fi
    \edef#1{\appendixalpha}
    \edef\currentsec{\appendixalpha}}

  \newcount\procno \procno=0
  \def\newproc#1{\advance\procno by 1
   \ifnum\ignorsec=0 \edef#1{\currentsec.\number\procno}
   \else \edef#1{\number\procno}
   \fi}

  \newcount\tagno \tagno=0
  \def\newtag#1{\advance\tagno by 1 \edef#1{\number\tagno}}



\notasec
  \newtag{\CstarFP}
\newsec{\MainRes}
  \newtag{\Hilsum}
 \newproc{\MainThm}
 \newproc{\VstaraV}
  \newtag{\aoneam}
  \newtag{\uglyRHS}
 \newproc{\VAVeqA}
  \newtag{\phiVaV}
\newsec{\ExampleSec}
 \newproc{\ExNotCyclic}

\newbib{\VoiculescuZZSymmetries}
\newbib{\VDNbook}

\topmatter

  \title Faithfulness of free product states
  \endtitle

  \leftheadtext{Ken Dykema}
  \rightheadtext{Faithfulness of free product states}

  \author Kenneth J\. Dykema
  \endauthor

  \date \today \enddate

  \affil 
    Department of Mathematics and Computer Science \\
    Odense Universitet, Campusvej 55 \\
    DK-5230 Odense M \\
    Denmark
  \endaffil

  \address
    Department of Mathematics and Computer Science,
    Odense Universitet, Campusvej 55,
    DK-5230 Odense M,
    Denmark
  \endaddress

  \abstract
    It is proved that the free product state, in the reduced free product of
    C$^*$--algebras, is faithful if the states one started with are faithful.
  \endabstract


\endtopmatter

\document \TagsOnRight \baselineskip=18pt

The reduced free product of C$^*$--algebras was introduced by
Voiculescu
as the appropriate construction for C$^*$--algebras in the setting of
his theory of freeness~\cite{\VoiculescuZZSymmetries}.
(See also the book~\cite{\VDNbook}).
Given a set $I$ and  for each $\iota\in I$ a unital C$^*$--algebra $A_\iota$
with a state $\phi_\iota$ whose GNS representation is faithful,
the reduced free product construction yields the unique, unital C$^*$--algebra
$A$ with unital embeddings $A_\iota\hookrightarrow A$ and a state $\phi$ on
$A$ such that
\roster
\item"(i)" $\phi\restrict_{A_\iota}=\phi_\iota$,
\item"(ii)" $(A_\iota)_{\iota\in I}$ is free in $(A,\phi)$,
\item"(iii)" $A=C^*(\bigcup_{\iota\in I}A_\iota)$,
\item"(iv)" the GNS representation of $\phi$ is faithful on $A$.
\endroster
We denote the reduced free product by
$$ (A,\phi)=\freeprodi(A_\iota,\phi_\iota) \tag{\CstarFP} $$
and $\phi$ is called the {\it free product state}.
In this note, we prove that if $\phi_\iota$ is faithful on $A_\iota$ for every
$\iota\in I$ then $\phi$ is faithful on $A$.

Voiculescu~\cite{\VoiculescuZZSymmetries} proved that the free product state in
the analogous construction for von~Neumann algebras is faithful
under the hypothesis that each $\phi_\iota$ is faithful.
This implies that for the reduced free product of C$^*$--algebras
in~(\CstarFP),
$\phi$ is faithful if, for each $\iota\in I$, letting
$(\pi_\iota,\Hil_\iota,\xi_\iota)$ denote the GNS construction for
$(A_\iota,\phi_\iota)$, the vector $\xi_\iota$ is cyclic for the commutant,
$\pi_\iota(A_\iota)'$ on $\Hil_\iota$.
This is always the case if $\phi_\iota$ is faithful and is a trace.
However, there
are $(A_\iota,\phi_\iota)$ with $\phi_\iota$ faithful but $\xi_\iota$ not
cyclic for $\pi_\iota(A_\iota)'$.
See example~\ExNotCyclic.

\vskip3ex
\noindent{\bf\S\MainRes.  The main result.}
\vskip3ex

Here is the main result, stated a little more formally.
\proclaim{Theorem \MainThm}
Let $I$ be a set and for each $\iota\in I$ let $A_\iota\ne\Cpx$ be a
unital C$^*$--algebra and let $\phi_\iota$ be a faithful state on $A_\iota$.
Let
$$ (A,\phi)=\freeprodi(A_\iota,\phi_\iota) $$
be the reduced free product of C$^*$--algebras.
Then the state $\phi$ is faithful on $A$.
\endproclaim
The hypothesis that $A_\iota\ne\Cpx$ is not really a restriction, because for
any C$^*$--algebra $B$ with state $\phi$, $(B,\phi)*(\Cpx,\id)=(B,\phi)$.

The proof of the theorem will rely on Voiculescu's construction of the reduced
free product,
which we now briefly review.
As above, let $(\pi_\iota,\Hil_\iota,\xi_\iota)=\GNS(A_\iota,\phi_\iota)$ and
let
$$ \Hilo_\iota=\Hil_\iota\ominus\Cpx\xi_\iota. $$
For $\zeta\in\Hil_\iota$ and $a\in A_\iota$ we will always write $a\zeta$
instead of $\pi_\iota(a)\zeta$.
Set
$$ \Hil=\Cpx\xi\oplus\bigoplus
\Sb n\in\Naturals\\\iota_1,\ldots,\iota_n\in I\\\iota_j\ne\iota_{j+1}\endSb
\Hilo_{\iota_1}\otimes\Hilo_{\iota_2}\otimes\cdots\otimes\Hilo_{\iota_n}.
\tag{\Hilsum} $$
In Voiculescu's construction, each $A_\iota$ is represented on $\Hil$ and $A$
is the C$^*$--algebra generated by the union these $A_\iota$ acting on $\Hil$.
Then $\phi$ is the vector state for $\xi$.

Let $n\in\Naturals$ and $\iota_1,\ldots,\iota_n\in I$ be such that
$\iota_1\ne\iota_2,\,\iota_2\ne\iota_3,\,\ldots,\iota_{n-1}\ne\iota_n$.
For each $1\le j\le n-1$ let $\zeta_j\in\Hilo_{\iota_j}$.
Define the isometry
$$ V=V_{(\zeta_1,\ldots,\zeta_{n-1},\iota_n)}:\Hil_{\iota_n}\to\Hil $$
by
$$ \align
V\xi_{\iota_n}&=\zeta_1\otimes\zeta_2\otimes\cdots\otimes\zeta_{n-1} \\
V\zeta&=\zeta_1\otimes\zeta_2\otimes\cdots\otimes\zeta_{n-1}\otimes\zeta
\quad(\zeta\in\Hilo_{\iota_n}).
\endalign $$

Whenever $B$ is a C$^*$--algebra and $\psi$ a state on $B$ we denote by $\bh$
the corresponding element in $L^2(B,\psi)$.
Also, we employ the standard notation $A_\iota\oup=\ker\phi_\iota$.
\proclaim{Lemma \VstaraV}
Let $V=V_{(\zeta_1,\ldots,\zeta_{n-1},\iota_n)}$ be as above.
Let $m\in\Naturals$, $k_1,\ldots,k_m\in I$ be such that
$k_1\ne k_2,\quad k_2\ne k_3,\quad\ldots,\quad k_{m-1}\ne k_m$ and let
$a_j\in A_{k_j}\oup$,
$(1\le j\le m)$.

If $m=2p-1$ for an integer $p$ with $1\le p<n$ and if
$$ k_m=\iota_1=k_1,\,k_{m-1}=\iota_2=k_2,\,\ldots,\,
k_{p+1}=\iota_{p-1}=k_{p-1},\,k_p=\iota_p $$
then $V^*a_1a_2\ldots a_mV=c1$ where $c$ is the scalar
$$ \align
c=&
\bigl(\langle a_m\zeta_1,\xi_{\iota_1}\rangle\,
\langle a_{m-1}\zeta_2,\xi_{\iota_2}\rangle\,\cdots\,
\langle a_{p+1}\zeta_{p-1},\xi_{\iota_{p-1}}\rangle\bigr)
\langle a_p\zeta_p,\zeta_p\rangle\;\cdot \\
&\cdot\;\bigl(\langle a_{p-1}\xi_{\iota_{p-1}},\zeta_{p-1}\rangle\,
\langle a_{p-2}\xi_{\iota_{p-2}},\zeta_{p-2}\rangle\,\cdots\,
\langle a_1\xi_{\iota_1},\zeta_1\rangle\bigr).
\endalign $$

If $m=2n-1$ and if
$$ k_m=\iota_1=k_1,\,k_{m-1}=\iota_2=k_2,\,\ldots,\,
k_{n+1}=\iota_{n-1}=k_{n-1},\,k_n=\iota_n $$
then $V^*a_1a_2\ldots a_mV=ca_n$ where $c$ is the scalar
$$ \aligned
c=&
\bigl(\langle a_m\zeta_1,\xi_{\iota_1}\rangle\,
\langle a_{m-1}\zeta_2,\xi_{\iota_2}\rangle\,\cdots\,
\langle a_{n+1}\zeta_{n-1},\xi_{\iota_{n-1}}\rangle\bigr)\;\cdot\\
&\cdot\;\bigl(\langle a_{n-1}\xi_{\iota_{n-1}},\zeta_{n-1}\rangle\,
\langle a_{n-2}\xi_{\iota_{n-2}},\zeta_{n-2}\rangle\,\cdots\,
\langle a_1\xi_{\iota_1},\zeta_1\rangle\bigr).
\endaligned \tag{\aoneam} $$

Otherwise, $V^*a_1a_2\ldots a_mV=0$.
\endproclaim
\demo{Proof}
We use induction on $n$.
If $n=1$ then $V(\Hil_{\iota_1})$ is the canonical copy
$\Hil_{\iota_1}\cong\Cpx\xi\oplus\Hilo_{\iota_1}$ of $\Hil_{\iota_1}$ in
$\Hil$.
If $m=1$ and $k_m=\iota_1$ then $V^*a_1V=a_1$, whereas if $m>1$ or if $m=1$ but
$k_m\ne1$ then one easily sees that
$$ \align
a_1\cdots a_m\xi=&\ah_1\otimes\cdots\otimes\ah_m
\quad\perp(\Cpx\xi\oplus\Hilo_{\iota_1}) \\
a_1\cdots a_m\Hilo_{\iota_1}=&\ah_1\otimes\cdots\otimes\ah_m
\otimes\Hilo_{\iota_1}
\quad\perp(\Cpx\xi\oplus\Hilo_{\iota_1}),
\endalign $$
so that $V^*a_1\cdots a_mV=0$.

Now suppose $n>1$.
If $k_m\ne\iota_1$ then
$$ \align
a_1\cdots a_m(\zeta_1\otimes\cdots\otimes\zeta_{n-1})
=&\ah_1\otimes\cdots\ah_m\otimes\zeta_1\otimes\cdots\otimes\zeta_{n-1}
\quad\perp V\Hil_{\iota_n} \\
a_1\cdots a_m(\zeta_1\otimes\cdots\otimes\zeta_{n-1}\otimes\Hilo_{\iota_n})
=&\ah_1\otimes\cdots\ah_m\otimes\zeta_1\otimes\cdots\otimes\zeta_{n-1}
\otimes\Hilo_{\iota_n}
\quad\perp V\Hil_{\iota_n}.
\endalign $$
So $k_m\ne\iota_1$ implies $V^*a_1\cdots a_mV=0$.
Taking adjoints we see that also $k_1\ne\iota_1$ implies
$V^*a_1\cdots a_mV=0$.
Hence suppose $k_m=\iota_1=k_1$.
Then $m=1$ or $m\ge3$.
Furthermore,
$$ \aligned
a_1\cdots a_m(\zeta_1\otimes\cdots\otimes\zeta_{n-1})
=&\ah_1\otimes\cdots\ah_{m-1}\otimes
\bigl(a_m\zeta_1-\langle a_m\zeta_1,\xi_{\iota_1}\rangle\xi_{\iota_1}\bigr)
\otimes\zeta_2\cdots\otimes\zeta_{n-1}+ \\
&+\langle a_m\zeta_1,\xi_{\iota_1}\rangle
a_1\cdots a_{m-1}(\zeta_2\otimes\cdots\otimes\zeta_{n-1})
\endaligned \tag{\uglyRHS} $$
and similarly for every element of
$\zeta_1\otimes\cdots\otimes\cdots\zeta_{n-1}\otimes\Hilo_{\iota_n}$.
If $m=1$ then the second term of the right hand side of~(\uglyRHS) is
perpendicular to $V\Hil_{\iota_n}$ and taking $V^*$ of~(\uglyRHS) shows
$$ V^*a_1V=\langle a_1\zeta_1,\zeta_1\rangle1. $$
If $m\ge3$ then the first term of the right hand side of~(\uglyRHS) is
perpendicular to
$V\Hil_{\iota_n}$ and we see that
$$ V^*a_1\cdots a_mV=\langle a_m\zeta_1,\xi_{\iota_1}\rangle
V^*a_1\cdots a_{m-1}U, $$
where $U=V_{(\zeta_2,\ldots,\zeta_{n-1},\iota_n)}$.
In a like manner we show that
$$ \align
V^*a_1\cdots a_{m-1}U&=(U^*a_{m-1}^*\cdots a_1^*V)^*
=\overline{\langle a_1^*\zeta_1,\xi_{\iota_1}\rangle}
(U^*a_{m-1}^*\cdots a_2^*U)^* \\
&=\langle a_1\xi_{\iota_1},\zeta_1\rangle U^*a_2\cdots a_{m-1}U.
\endalign $$
Now applying the inductive hypothesis finishes the proof.
\QED

\proclaim{Lemma \VAVeqA}
Let $(A,\phi)$ be as in Theorem~\MainThm.
Let $n\in\Naturals$ and $\iota_1,\ldots,\iota_n\in I$ be such that
$\iota_1\ne\iota_2,\,\iota_2\ne\iota_3,\,\ldots,\iota_{n-1}\ne\iota_n$.
For each $1\le j\le n-1$ let $\zeta_j\in\Hilo_{\iota_j}$.
Let $V=V_{(\zeta_1,\ldots,\zeta_{n-1},\iota_n)}$ be as defined before
Lemma~\VstaraV.
Then $V^*AV=A_{\iota_n}$.
\endproclaim
\demo{Proof}
Clearly $V^*V=1\in A_{\iota_n}$ and from Lemma~\VstaraV,
$V^*a_1\cdots a_mV\in A_{\iota_n}$ for each
$a_1\cdots a_m$ as in the statement of Lemma~\VstaraV.
Since the collection of all these $a_1\cdots a_m$ together with $1$ spans a
dense subspace of $A$, it follows that $V^*AV\subseteq A_{\iota_n}$.
In order to see the other inclusion, using~(\aoneam) it is enough to see that
for every $1\le j\le n-1$ one can find $a_j,a_{2n-j}\in A_{\iota_j}\oup$ such
that
$$ \langle a_j\xi_{\iota_j},\zeta_j\rangle\ne0,\qquad
\langle a_{2n-j}\zeta_j,\xi_{\iota_j}\rangle\ne0. $$
But this can be done because $\xi_{\iota_j}$ is cyclic for the action of
$A_{\iota_j}$ on $\Hil_{\iota_j}$.
\QED

\demo{Proof of Theorem~\MainThm}
Suppose for contradiction there is $a\in A$ with $a\ge0$, $a\ne0$ and
$\phi(a)=0$.
Then
$$ \langle a\xi,\xi\rangle=\phi(a)=0. $$
If $n\in\Naturals$ and $\iota_1,\ldots,\iota_n\in I$ with
$\iota_1\ne\iota_2,\,\iota_2\ne\iota_3,\,\ldots,\iota_{n-1}\ne\iota_n$, then
let $p_{\iota_1,\iota_2,\ldots,\iota_n}$ denote the projection from $\Hil$ onto
the direct summand
$$ \Hilo_{\iota_1}\otimes\Hilo_{\iota_2}\otimes\cdots\otimes\Hilo_{\iota_n} $$
in~(\Hilsum).
Since $a\ge0$ and $a\ne0$, for some $n\in\Naturals$ one can choose
$\iota_1,\ldots,\iota_n$ as above such that
$$ p_{\iota_1,\ldots,\iota_n}ap_{\iota_1,\ldots,\iota_n}\ne0. $$
Let $n$ be least for which such a choice can be made.
Then there are $\zeta_j\in\Hilo_{\iota_j}$, $(1\le j\le n-1)$, such that
letting $V=V_{(\zeta_1,\ldots,\zeta_{n-1},\iota_n)}$ we have $V^*aV\ne0$.
However, $V^*aV\ge0$ and by Lemma~\VAVeqA{} $V^*aV\in A_{\iota_n}$.
Since $\phi_{\iota_n}$ is faithful on $A_{\iota_n}$, we must have
$$ 0<\phi_{\iota_n}(V^*aV)=\langle V^*aV\xi_{\iota_n},\xi_{\iota_n}\rangle.
\tag{\phiVaV} $$
If $n=1$ then
$$ \langle V^*aV\xi_{\iota_n},\xi_{\iota_n}\rangle
=\langle a\xi,\xi\rangle=\phi(a)=0, $$
contradicting~(\phiVaV).
If $n>1$ then from~(\phiVaV) and the definition of $V$ we have
$$ 0<\langle V^*aV\xi_{\iota_n},\xi_{\iota_n}\rangle
=\langle a(\zeta_1\otimes\cdots\otimes\zeta_{n-1}),
\zeta_1\otimes\cdots\otimes\zeta_{n-1}\rangle. $$
Hence $p_{\iota_1,\ldots,\iota_{n-1}}ap_{\iota_1,\ldots,\iota_{n-1}}\ne0$,
which contradicts the choice of $n$.
\QED

\vskip3ex
\noindent{\bf\S\ExampleSec.  An example.}
\vskip3ex

In this section we present an example of a C$^*$--algebra $A$ with a faithful
state $\phi$ such that in the GNS representation the standard vector is not
cyclic for the action of the commutant.

\proclaim{Example \ExNotCyclic}\rm
Let $\Teu$ be the Toeplitz algebra generated by the nonunitary isometry, $S$,
and let $\pi:\Teu\to C(\Tcirc)$ be a unital $*$--homomorphism.
Let $\ev:C(\Tcirc)\to\Cpx$ be a character of $C(\Tcirc)$ and let
$\psi_0=\ev\circ\pi$.
Let $\psi_1$ be any faithful state on $\Teu$.
Let $\tr_2$ be the tracial state on $M_2(\Cpx)$ and let $\rho$ be a pure state
on $M_2(\Cpx)$.
Let $A=\Teu\otimes M_2(\Cpx)$ and let
$$ \phi=\tfrac12\psi_1\otimes\tr_2+\tfrac12\psi_0\otimes\rho:A\to\Cpx. $$
Clearly, $\phi$ is a faithful state on $A$.
However, letting
$$ (\pi,\Hil,\xi)=\GNS(A,\phi), $$
we will show that $\pi(A)'\xi\ne\Hil$.

Our first task is to find $\Hil$.
Let $(e_{ij})_{i,j\in\Naturals}$ be a system of matrix units for the copy,
$\KHil$, of the compact operators which is an ideal of $\Teu$, and let
$(f_{ij})_{1\le i,j\le2}$ be a system of matrix units for $M_2(\Cpx)$ such that
$\rho(f_{11})=1$.
Then
$$ \align
L^2(A,\psi_1\otimes\tr_2)\text{ has orthonormal basis }
&\{\psi_1(e_{jj})^{-1/2}\eh_{ij}\mid i,j\in\Naturals\}\otimes
\{\sqrt2\fh_{ij}\mid1\le i,j\le2\} \\
L^2(A,\psi_0\otimes\rho)\text{ has orthonormal basis }
&\onehat\otimes\{\fh_{11},\fh_{21}\}.
\endalign $$
Let
$$ \Hil'=L^2(A,\psi_1\otimes\tr_2)\oplus L^2(A,\psi_0\otimes\rho) $$
with norm
$\nm{\zeta_1\oplus\zeta_0}^2=\frac12(\nm{\zeta_1}^2+\nm{\zeta_0}^2)$.
Since for every $a\in A$ the elements $\ah\in L^2(A,\phi)$ and
$\ah\oplus\ah\in\Hil'$ have the
same norm, there is an isometry $V:L^2(A,\phi)\to\Hil'$
such that $V\ah=\ah\oplus\ah$.
Let us show that $V$ is onto $\Hil'$.
Clearly
$$ V(e_{ij}\otimes f_{kl})\hat{\;}=(\eh_{ij}\otimes\fh_{kl})\oplus0. $$
Thus $L^2(A,\psi_1\otimes\tr_2)\oplus0$ is in the range of $V$.
But for $i=1,2$
$$ V(\onehat\otimes\fh_{i1})\hat{\;}=(\onehat\otimes\fh_{i1})\oplus
(\onehat\otimes\fh_{i1}), $$
so $0\oplus(\onehat\otimes\fh_{i1})$ is in the range of $V$.
Hence $V$ is onto $\Hil'$.
Consequently, we may identify $\Hil$ with $\Hil'$ and $\xi\in\Hil$ with
$$ (\onehat\otimes\onehat)\oplus(\onehat\otimes\eh_{ii})\in\Hil'. $$
We see that $L^2(\Teu,\psi_1)$ can be identified with
$l^2(\Naturals)\otimes l^2(\Naturals)$ by
$$ L^2(\Teu,\psi_1)\ni\psi_1(e_{jj})^{-1/2}\eh_{ij}
\mapsto \delta_i\otimes\delta_j, $$
and the left action of $\Teu$ on $L^2(\Teu,\psi_1)$ is realized on
$l^2(\Naturals)\otimes l^2(\Naturals)$ by
$S(\delta_i\otimes\delta_j)=\delta_{i+1}\otimes\delta_j$.
Moreover, identifying $L^2(M_2(\Cpx),\tr_2)$ with
$l^2(\{1,2\})\otimes l^2(\{1,2\})$ by
$\sqrt2\fh_{ij}\mapsto\delta_i\otimes\delta_j$,
and identifying $L^2(A,\psi_0\otimes\rho)$ with $l^2(\{1,2\})$ in the obvious
way, we identify $\Hil$ with
$$ \bigl(l^2(\Naturals)\otimes l^2(\{1,2\})\otimes l^2(\{1,2\}) \otimes
l^2(\Naturals)\bigr)\oplus l^2(\{1,2\}). $$

Now we find $\pi(A)'$ by first finding $\pi(A)''$.
Clearly
$$ \pi(A)''\supseteq\pi(\KHil\otimes M_2(\Cpx))''
=\bigl(B(l^2(\Naturals))\otimes B(l^2(\{1,2\}))\otimes1\otimes1\bigr)\oplus0 $$
Moreover, considering now $\pi(1\otimes M_2(\Cpx))$ we see that
$$ \pi(A)''=
\bigl(B(l^2(\Naturals))\otimes B(l^2(\{1,2\}))\otimes1\otimes1\bigr)
\oplus B(l^2(\{1,2\})). $$
Therefore
$$ \pi(A)'=
\bigl(1\otimes1\otimes B(l^2(\{1,2\}))\otimes B(l^2(\Naturals))\bigr)
\oplus\Cpx $$
and thus
$0\oplus(\onehat\otimes\fh_{21})\perp\pi(A)'\xi$,
which was what we wanted to show.
\endproclaim

\enddocument